\documentclass{article}

\usepackage{arxiv}

\usepackage[utf8]{inputenc} % allow utf-8 input
\usepackage[T1]{fontenc}    % use 8-bit T1 fonts
\usepackage{hyperref}       % hyperlinks
\usepackage{url}            % simple URL typesetting
\usepackage{booktabs}       % professional-quality tables
\usepackage{amsfonts}       % blackboard math symbols
\usepackage{nicefrac}       % compact symbols for 1/2, etc.
\usepackage{microtype}      % microtypography
\usepackage{lipsum}
\usepackage{graphicx}
\graphicspath{ {./images/} }
\usepackage{amsmath}

\title{Characterization of\(\left(\alpha,\alpha\right)\)-derivation on B(X)}

\author{
 Quanyuan Chen \\
  School of Information and Engineering,\\
  Jingdezhen Ceramic University, \\
  Jingdezhen 333000. China. \\
  \texttt{cqy0798@163.com} \\
   \And
 Yaqi Li \\
  School of Information and Engineering,\\
 Jingdezhen Ceramic University,\\
 Jingdezhen 333000. China. \\
  \texttt{lyq8532@163.com} \\
}

\begin{document}
\maketitle
\begin{abstract}
Let \( X\) be a Banach algebra and \( B(X)\) be the set of all bounded linear operators on \( X\). Suppose that \(\alpha: B(X)\rightarrow B(X)\) is an automorphism. We say that a mapping \(\delta\) from \(B(X)\) into itself is derivable at \(G\in B(X)\) if \(\delta(G) =\alpha\left(A\right)\delta(B) + \delta(A)\alpha\left (B\right)\) for all \(A\), \(B\in B(X)\) with \(AB = G\). We say that an element \(G\in B(X)\) is an \(\left(\alpha,\alpha\right)\)-all derivable point of \(B(X)\) if every \(\left(\alpha,\alpha\right)\)-derivable mapping \(\delta\) at \(G\) is an \(\left(\alpha,\alpha\right)\)- derivation. In this paper, we show that every \(\left(\alpha,\alpha\right)\)-derivable mapping at a nonzero element in \(B(X)\) is an \(\left(\alpha,\alpha\right)\)-derivation. 
\end{abstract}

% keywords can be removed

\section{Introduction}
Let  \(\mathcal{A}\)  be an algebra over  \(\mathbb{F}\) , where  \(\mathbb{F}\)  is either the real field  \(\mathbb{R}\)  or the complex field  \(\mathbb{C}\) , and  \(\delta: \mathcal{A} \rightarrow \mathcal{A}\)  be an additive (linear) mapping. Recall that  \(\delta\)  is called a derivation if  \(\delta(A B)=\delta(A) B+A \delta(B)\)  for all  \(A, B \in \mathcal{A}\) , and  \(\delta\)  is called a Jordan derivation if  \(\delta\left(A^{2}\right)=\delta(A) A+A \delta(A)\)  for all  \(A \in \mathcal{A}\) . Suppose that  \(\alpha\)  and  \(\beta\)  are ring automorphisms on  \(\mathcal{A}\). An additive mapping  \(\delta: \mathcal{A} \rightarrow \mathcal{A}\)  is called an additive  \((\alpha, \beta)\) -derivation if  \(\delta(A B)=\alpha(A) \delta(B)+\delta(A) \beta(B)\)  for all  \(A, B \in \mathcal{A}\) . An additive mapping  \(\delta: \mathcal{A} \rightarrow \mathcal{A}\)  is called an additive  \((\alpha, \beta)\)  - generalized derivation if  \(\delta(A B)=\alpha(A) \delta(B)+\delta(A) \beta(B)-\alpha(A) \delta(I) \beta(B)\)  for all  \(A, B \in \mathcal{A}\) .

Derivations are very important mappings both in theory and applications and have been studied extensively. In \cite{jacobson1950jordan}, Jacobson and Rickart proved that every Jordan derivation on an unital ring  \(\mathrm{R}\)  is a derivation. In \cite{r2}, Herstein proved that every Jordan derivation on a 2-torsion free prime ring is a derivation. In \cite{r3}, Xue et al. proved that if \(\mathcal{A}=Alg\mathcal{N}\) is a nest algebra and  \(\mathcal{M}=B(X)\), then a linear mapping \(\delta: \mathcal{A} \rightarrow \mathcal{M}\) is Jordan derivable at a nontrivial projection  \(P \in Alg\mathcal{N}\) if and only if it is derivation. In \cite{r4}, Jing and Lu studied \(\left(\alpha, \beta\right)\)-derivations in operator algebras and characterized \(\left(\alpha, \beta\right)\)-derivations on \(B(X)\). 

In recent years, some scholars have been interested in characterizing the mapping as derivable at some point. Suppose that \(G\) is an element in \(\mathcal{A}\). We say that  \(\delta\) is a derivable mapping at \(G\) if \(\delta(G)=  \delta(A) B+A \delta(B)\) for all \(A, B \in \mathcal{A}\) with \(AB=G\). An element \(G \in \mathcal{A}\) is called an all-derivable point in  \(\mathcal{A}\) if every derivable mapping at \(G\) is a derivation. Jing et al. in \cite{r5} proved that every linear mapping derivable at zero point with \(\delta(I)=0\) is an inner derivation of the nest algebra on Hilbert spaces. Zhu proved in \cite{r6} that every invertible operator on the nest algebra \(Alg\mathcal{N}\) is an all-derivable point of the nest algebra under the strong operator topology. In \cite{r7}, Jing proved that the zero point is a generalized Jordan all-derivable point of \(B(H)\) if  H  is infinite-dimensional and  I  is a Jordan all-derivable point of \(B(H)\) of a Hilbert space. He et al. proved in \cite{r8} that every derivable mapping at a nonzero element is a derivation in \(B(H)\). 

Inspired by \cite{r8}, we consider \(\left(\alpha,\alpha\right)\)-all derivable points for an automorphism  \(\alpha\) on \(\mathcal{A}\) . We say that \(\delta\) is \(\left(\alpha,\alpha\right)\)-derivable at point \(G \in \mathcal{A}\) if \(\delta(G)=\alpha(A) \delta(B)+\delta(A) \alpha(B)\) for all \(A, B \in \mathcal{A}\) with \(AB=G\). It is obvious that \((I, I)\)-derivations are usual derivations, and  \(\left(\alpha,\alpha\right)\)-all derivable points are all derivable points if \(\alpha\) is the identity mapping. In \cite{r9}, Jiao and Hou showed that if \(\delta\) is generalized \((\alpha, \beta)\)-derivable at the zero, then \(\delta\) is an additive generalized \((\alpha, \beta)\)-derivation, where \(Alg\mathcal{N}\) is a nest algebra on a (real or complex) Banach space.

The purpose of this paper is to characterize the mapping \((\alpha, \alpha)\)-derivable at a nonzero element in \(B(X)\). In section 2, we prove that every \((\alpha, \alpha)\)-derivable mapping at a nonzero element in \(B(X)\) is an \((\alpha, \alpha)\)-derivation. Moreover, we also prove that every \((\alpha, \alpha)\)-generalized derivable mapping at a nonzero element \(G\) in \(B(X)\) (that is, \(\delta(G)=\delta(A) \alpha(B)+\alpha(A) \delta(B)-\alpha(A) \delta(I) \alpha(B)\) for all \(A, B \in \mathcal{A}\) with  \(A B=G\).) is a \((\alpha, \alpha)\)-general derivation.

The following lemma will get repeated use.

\textbf{Lemma 1} (Lemma 1.1 in \cite{r8}). Let \(\mathcal{R}\) be a semiprime ring. \(b\), \(p\), \(q\) are fixed elements in \(\mathcal{R}\). 

(1) If \(babab= 0\) for all \(a\in \mathcal{R}\), then \(b=0\).

(2) If \(paqbpaq= 0\) for all \(a\in \mathcal{R}\), then \(paq=0\). 

\textbf{Lemma 2} (Corollary 1 in \cite{r10}). Let \(\mathcal{R}\) be a 2-torsion free semiprime ring and let \(\alpha\), \(\beta\) be an automorphism of \(\mathcal{R}\). Then every Jordan \((\alpha, \beta)\)-derivation of \(\mathcal{R}\) is an \((\alpha, \beta)\)-derivation.

\section{Main theorem}

In this section, \(X\) is a Banach space on \(\mathbb{F}\), where \(\mathbb{F}\) is either the real field \(\mathbb{R}\) or the complex field \(\mathbb{C}\). \(B(X)\) is an algebra of all bounded linear operators on \(X\). For convenience, we use the symbol \(\mathcal{A}\) to represent \(B(X)\) in the following proof. \(I\) is the identity element of \(\mathcal{A}\), and \(X^{\ast}\) is the set of all bounded linear functionals from \(X\) into \(\mathbb{F}\).

For any \textit{}\(x\in X\), \(f\in X^{\ast}\), the rank one operator on \(X\) is defined as \(x\otimes f\left(y\right)=f\left(y\right)x\) for all \(y\in X\). Clearly, \(x\otimes f\) is a rank one idempotent operator when \(f\left(x\right)=1\). In this case, we have \(\left(x\otimes f\right)A\left(x\otimes f\right)=\lambda x\otimes f\) for some \(\lambda= f\left(AX\right)\).

According to the Peirce decomposition, for a nontrivial idempotent \(P_{1}\) in \(\mathcal{A}\) and \(P_{2}=I-P_{1}\), we denote \(\mathcal{A}_{11} = P_{1} \mathcal{A} P_{1}\), \(\mathcal{A}_{12} = P_{1} \mathcal{A} P_{2}\), \(\mathcal{A}_{21} = P_{2} \mathcal{A} P_{1}\), \(\mathcal{A}_{22} = P_{2} \mathcal{A} P_{2}\), respectively, then \(\mathcal{A}\) can be decomposed as \(\mathcal{A} = \mathcal{A}_{11}+\mathcal{A}_{12}+\mathcal{A}_{21}+\mathcal{A}_{22}\). For any operator \(A\in\mathcal{A}\), it can be decomposed as \(A=A_{11}+A_{12}+A_{21}+A_{22}\), where \(A_{ij}\in\mathcal{A}_{ij}\), \(i\),\( j= 1\),\( 2\). Inspired by [10], we consider \((\alpha, \alpha)\)-all derivable points for an automorphism  \(\alpha\) on \(\mathcal{A}\). Let \(Q_{1} = \alpha\left(P_{1}\right)\) and \(Q_{2} = \alpha\left(P_{2}\right)\) , then \(Q_{1}\) and \(Q_{2}\) are nontrivial idempotents in \(\mathcal{A}\) with \(Q_{2} = I-Q_{1}\).

The following are our main results of the paper.

\textbf{Theorem 2.1.} Suppose that \(G\) is a nonzero element in \(\mathcal{A} = B(X)\),  \(\alpha\) is an automorphism on \(\mathcal{A}\) , and \(\delta\) is a linear mapping   on \(\mathcal{A}\). If \(\delta\) is \((\alpha, \alpha)\)-derivable at \(G\), which  means that \(\delta(G) = \delta(A) \alpha(B)+\alpha(A) \delta(B)\) for all \(A, B\in \mathcal{A}\) with \(AB=G\) , then \(\delta\) is an \((\alpha, \alpha)\)-derivation.

\textbf{Proof.} If \(\overline{ran\alpha(G)} = X\), where \(\overline{ran\alpha(G)}\)   is the closure of the range of \(\alpha(G)\), then \(\alpha(G)\) is a right separating point of \(B(X)\). By \cite{r11}, \(\delta\) is a Jordan \((\alpha, \alpha)\)-derivation, then  \(\delta\) is an \((\alpha, \alpha)\)-derivation by Lemma 2.

Next, it is only necessary to prove the case in which \(\overline{ran\alpha(G)} \neq X\). Take a nonzero element \(\alpha\left(x_{0}\right)\in X\) which is not in \(\alpha(G)\). According to the Hahn-Banach Theorem, consider \(f_{0} \in \alpha\left(X^{\ast}\right)\) such that \(f_{0}\left(\alpha\left(x_{0}\right)\right) = 1\) and \(f_{0}(\operatorname{ran}  \alpha(G)) = 0\). Let \(Q_{2} = \alpha\left(x_{0}\right) \otimes f_{0}\), \(Q_{1} = I-Q_{2}\), then we have \(Q_{2} \alpha(G) = 0\), \(Q_{1} \alpha(G) = \alpha(G)\).

Assume that  \(A\), \(B\) are arbitrary elements in \(\mathcal{A}\), \(A_{ij}\), \(B_{ij}\) are any element in \(\mathcal{A}_{ij}\), \(i, j=1,2\), and \(t, s\) are any nonzero elements in \(\mathbb{F}\). According to the condition that \(Q_{1}\) is surjective on \(\alpha\left(\mathcal{A}_{11}\right)\), it is known that \(A_{11}\) is invertible in \(\mathcal{A}_{11}\), which provides the next decomposition of \(G\).

First, we define a new linear mapping  \(\phi: \mathcal{A} \rightarrow \mathcal{A}\)  by \(\phi(A)=\delta(A)+[\alpha(A), D]\), where \(D=\alpha\left(P_{1}\right) \delta\left(P_{2}\right) \alpha\left(P_{2}\right)-\alpha\left(P_{2}\right) \delta\left(P_{2}\right) \alpha\left(P_{1}\right)\). According to the assumption that  \(\delta: \mathcal{A} \rightarrow \mathcal{A}\) is \((\alpha, \alpha)\)-derivable at the point \(G\), \(\phi\) is also an \((\alpha, \alpha)\)-derivable at the point \(G\) by linear operations and satisfies \(\phi\left(P_{2}\right) \in \alpha\left(\mathcal{A}_{11}\right)+\alpha\left(\mathcal{A}_{22}\right)\). Now, it is only necessary to prove that \(\phi\) is an \((\alpha, \alpha)\)-derivation to obtain that \(\delta\) is an \((\alpha, \alpha)\)derivation.
The point \(G\) is written as:

\begin{equation*}
    (A_{11}+tsA_{11}A_{12})(A_{11}^{-1}G-sA_{12}A_{2i}+t^{-1}A_{2i})=G, i=1,2.
\end{equation*}

According to \(\phi\) is \((\alpha, \alpha)\)-derivable at the point \(G\) gives:
\begin{equation*}
\phi(G)=\phi\left(A_{11}+t s A_{11} A_{12}\right) \alpha\left(A_{11}^{-1} G-s A_{12} A_{2 i}+t^{-1} A_{2 i}\right)
\end{equation*}
\begin{equation}
+\alpha\left(A_{11}+t s A_{11} A_{12}\right) \phi\left(A_{11}^{-1} G-s A_{12}A_{2i}+t^{-1} A_{2 i}\right) .
\end{equation}

According to the linear property, we expand Equation 1 and simplify the coefficients of \(t^{-1}\), \(s\), and \(ts^{2}\), respectively. We obtain the following three equations:

\begin{equation}
\phi\left(A_{11}\right) \alpha\left(A_{2 i}\right)+\alpha\left(A_{11}\right) \phi\left(A_{2 i}\right)=0 
\end{equation}
\begin{equation}
\phi\left(A_{11} A_{12}\right) \alpha\left(A_{2 i}\right)+\alpha\left(A_{11} A_{12}\right) \phi\left(A_{2 i}\right)=\phi\left(A_{11}\right) \alpha\left(A_{12} A_{2 i}\right)+\alpha\left(A_{11}\right) \phi\left(A_{12} A_{2 i}\right) 
\end{equation}
\begin{equation}
\phi\left(A_{11} A_{12}\right) \alpha\left(A_{12} A_{2 i}\right)+\alpha\left(A_{11} A_{12}\right) \phi\left(A_{12} A_{2 i}\right)=0 .
\end{equation}

By the linear property of \(\phi\) and the decomposition of arbitrary element \(A\) and \(B\) in \(\mathcal{A}\), to prove that \(\phi\) is an \((\alpha, \alpha)\)-derivation it is sufficient to show that the following equation holds:
\begin{equation}
\phi\left(A_{i j} B_{s k}\right)=\phi\left(A_{i j}\right) \alpha\left(B_{s k}\right)+\alpha\left(A_{i j}\right) \phi\left(B_{s k}\right)
\end{equation}

where  \(1 \leq i, j, s, k \leq 2 \).

\textbf{Claim 1}: \( \phi\left(A_{22}\right)=0\) .

We consider  \(A_{11}=P_{1}\)  and  \(A_{2 i}=P_{2}\)  in Equation 2, and obtain  \(\phi\left(P_{1}\right) \alpha\left(P_{2}\right)+   \alpha\left(P_{1}\right) \phi\left(P_{2}\right)=0\). We multiply the right-hand side of the above equation by  \(Q_{1}=\alpha\left(P_{1}\right)\)  and obtain  \(\alpha\left(P_{1}\right) \phi\left(P_{2}\right) \\\alpha\left(P_{1}\right)=0 \). By the known condition  \(\phi\left(P_{2}\right) \in \alpha\left(\mathcal{A}_{11}\right)+\alpha\left(\mathcal{A}_{22}\right)\) , we obtain  \(\phi\left(P_{2}\right) \in \alpha\left(\mathcal{A}_{22}\right) \), so  \(\phi\left(P_{2}\right)=\lambda \alpha\left(P_{2}\right) \) for some  \(\lambda \in \mathbb{F} \) and  \(\phi\left(P_{1}\right) \alpha\left(P_{2}\right)=0 \). We consider  \(A_{2 i}=P_{2} \) in Equation 2 again and obtain
\begin{equation}
\phi\left(A_{11}\right) \alpha\left(P_{2}\right)=0 .
\end{equation}
Letting  \(A_{11}=P_{1}\)  and  \(A_{2 i}=A_{21}\)  in Equation 2, we obtain that
\begin{equation}
\alpha\left(P_{1}\right) \phi\left(A_{21}\right)=0 \text {. }
\end{equation}

Cunsidering that  \(A_{11}=P_{1}\)  and \( A_{2 i}=P_{2} \) in Equation 3 and  \(\phi\left(P_{2}\right)=\lambda \alpha\left(P_{2}\right)\), \(\lambda \in \mathbb{F}\) , we get that
\begin{equation}
\phi\left(A_{12}\right) \alpha\left(P_{2}\right)+\lambda \alpha\left(A_{12}\right)=\phi\left(P_{1}\right) \alpha\left(A_{12}\right)+\alpha\left(P_{1}\right) \phi\left(A_{12}\right) .
\end{equation}

Multiply the right-hand side of Equation 8 by \( Q_{1}=\alpha\left(P_{1}\right) \), we obtain that
\begin{equation}
\alpha\left(P_{1}\right) \phi\left(A_{21}\right) \alpha\left(P_{1}\right)=0 \text {. }
\end{equation}

We multiply the left-hand side of Equation 8 by \( Q_{1}=\alpha\left(P_{1}\right)\) , and associated with Equation 9, we obtain \( \alpha\left(P_{1}\right) \phi\left(P_{1}\right) \alpha\left(A_{12}\right)=\lambda \alpha\left(A_{12}\right)\) . It is known that  \(\mathcal{A} \) is simple, so the  \(\alpha(\mathcal{A}) \) is also simple, and we get:
\begin{equation}
\alpha\left(P_{1}\right) \phi\left(P_{1}\right) \alpha\left(P_{1}\right)=\lambda \alpha\left(P_{1}\right) .
\end{equation}

As \( I G=G \) and \( \phi \) is an \( (\alpha, \alpha)\) -derivable at point  G , we have  \(\phi(I) \alpha(G)+\alpha(I) \phi(G)=   \phi(G) \), then \( \phi(I) \alpha(G)=0 \). Since  \(Q_{2} \alpha(G)=\alpha\left(P_{2}\right) \alpha(G)=0\)  and  \(\phi\left(P_{2}\right)=\lambda \alpha\left(P_{2}\right) \), then we have  \(\phi\left(P_{1}\right) \alpha(G)=\phi(I) \alpha(G)-\phi\left(P_{2}\right) \alpha(G)=-\lambda \alpha\left(P_{2}\right) \alpha(G)=0\) . We multiply the righthand side of Equation 10 by  \(\alpha(G)\) , due to \( Q_{1} \alpha(G)=\alpha\left(P_{1}\right) \alpha(G)=\alpha(G)\) , we get  \(\lambda \alpha(G)=\lambda \alpha (p_{1})\alpha(G)= \lambda \alpha (P_{1})\phi (P_{1})\alpha(G)= 0\). Since \(G\neq0\) and \(\alpha(G)\neq0\), which in turn gives \(\lambda=0\). Now we prove:
\begin{equation}
\alpha\left(P_{1}\right)\phi\left (P_{1}\right)\alpha\left (P_{1} \right)=\lambda \alpha\left(P_{1}\right)=0
\end{equation}
and \(\phi\left(P_{2}\right)=\lambda \alpha(P_{2})=0\), so we obtain \(\phi(A_{22})=0\).

\textbf{Claim 2}:\(\phi(\mathcal{A}_{11})\subseteq \alpha(\mathcal{A}_{11})\).

We consider \({A}_{11}=P_{1}\)  and \({A}_{2i}=P_{2}\)  in Equation 4, and we get
\begin{equation}
\phi\left(A_{12}\right) \alpha\left(A_{12}\right)+\alpha\left(A_{12}\right) \phi\left(A_{12}\right)=0 .
\end{equation}
The left multiplication of Equation 12 by \( Q_{1}=\alpha\left(P_{1}\right) \) and associated with Equation 9, then obtain:
\begin{equation}
\alpha\left(A_{12}\right) \phi\left(A_{12}\right)=0
\end{equation}
and
\begin{equation}
\phi\left(A_{12}\right) \alpha\left(A_{12}\right)=0 \text {. }
\end{equation}

According to Equation 13, we can obtain \( \alpha\left(B_{12}\right) \phi\left(B_{12}\right)=0  \) and \( \alpha\left(A_{12}+B_{12}\right) \phi\left(A_{12}+\right.   \left.B_{12}\right)=0 \), then we have:
\begin{equation}
\alpha\left(A_{12}\right) \phi\left(B_{12}\right)+\alpha\left(B_{12}\right) \phi\left(A_{12}\right)=0 .
\end{equation}

We multiply the right-hand side of Equation 15 by  \(\alpha\left(B_{12}\right) \) and associated with Equation 14 , we can obtain \( \alpha\left(B_{12}\right) \phi\left(A_{12}\right) \alpha\left(B_{12}\right)=0 \). By Lemma 1 we have  \(\alpha\left(P_{2}\right) \phi\left(A_{12}\right) \alpha\left(P_{1}\right)=0\) . Associating with Equation 9, we get:
\begin{equation}
\phi\left(A_{12}\right) \alpha\left(P_{1}\right)=0 .
\end{equation}

We consider \( A_{2 i}=P_{2} \) in Equation 3 then associated with  \(\phi\left(P_{2}\right)=0 \) and Equation 16 yields:
\begin{equation}
\phi\left(A_{11} A_{12}\right)=\phi\left(A_{11}\right) \alpha\left(A_{12}\right)+\alpha\left(A_{11}\right) \phi\left(A_{12}\right) .
\end{equation}

So, we have:

\begin{equation*}
\begin{aligned}
\phi\left(A_{11} B_{11} A_{12}\right) & =\phi\left(A_{11} B_{11}\right) \alpha\left(A_{12}\right)+\alpha\left(A_{11} B_{11}\right) \phi\left(A_{12}\right) \\
&=\phi\left(A_{11}\right) \alpha\left(B_{11} A_{12}\right)+\alpha\left(A_{11}\right) \phi\left(B_{11} A_{12}\right) \\
&=\phi\left(A_{11}\right) \alpha\left(B_{11} A_{12}\right)+\alpha\left(A_{11}\right) \phi\left(B_{11}\right) \alpha\left(A_{12}\right)+\alpha\left(A_{11} B_{11}\right) \phi\left(A_{12}\right).
\end{aligned}
\end{equation*}
Since  \(\phi\left(A_{11}\right) \alpha\left(P_{2}\right)=0\) , the above equation implies that
\begin{equation}
\phi\left(A_{11} B_{11}\right)=\phi\left(A_{11}\right) \alpha\left(B_{11}\right)+\alpha\left(A_{11}\right) \phi\left(B_{11}\right) .
\end{equation}
Letting \( A_{11}=P_{1} \) in Equation 18 we obtain that
\begin{equation}
\phi\left(B_{11}\right)=\phi\left(P_{1}\right) \alpha\left(B_{11}\right)+\alpha\left(P_{1}\right) \phi\left(B_{11}\right) .
\end{equation}

We consider \( A_{11}=P_{1} \) and \( A_{2 i}=A_{21} \) in Equation 4, associated with Equation 16 we obtain \( \alpha\left(A_{12}\right) \phi\left(A_{12} A_{21}\right)=0 \) and associated with Equation 19 we obtain \( \alpha\left(A_{12}\right) \phi\left(P_{1}\right) \alpha\left(A_{12} A_{21}\right)=0\). Since \( \alpha(\mathcal{A}) \) is simple, we have  \(\alpha\left(A_{12}\right) \phi\left(P_{1}\right) \alpha\left(A_{12}\right)=0 \). By Lemma 1 and  \(\alpha\)  is an automorphism on \( \mathcal{A}\) we have that \( \alpha\left(P_{2}\right) \phi\left(P_{1}\right) \alpha\left(P_{1}\right)=0 \). We multiply the left-hand side of Equation 19 by \( Q_{2}=\alpha\left(P_{2}\right)\)  obtains \( \alpha\left(P_{2}\right) \phi\left(B_{11}\right)=0 \). Together with \( \phi\left(B_{11}\right) \alpha\left(P_{2}\right)=0 \), this implies that  \(\phi\left(\mathcal{A}_{11}\right) \subseteq \alpha\left(\mathcal{A}_{11}\right) \).

\textbf{Claim 3}:  \(\phi\left(\mathcal{A}_{12}\right) \subseteq \alpha\left(\mathcal{A}_{12}\right)\) .

We consider  \(A_{11}=P_{1} \) in Equation 17 and multiply the left-hand side by \( Q_{2}=\alpha\left(P_{2}\right)\) , associated with  \(\phi\left(\mathcal{A}_{11}\right) \subseteq \alpha\left(\mathcal{A}_{11}\right) \), then we wet that \( \alpha\left(P_{2}\right) \phi\left(A_{12}\right)=0 \). Toghter with Equation 6:  \(\phi\left(A_{12}\right) \alpha\left(P_{1}\right)=0 \), we obtain that  \(\phi\left(\mathcal{A}_{12}\right) \subseteq \alpha\left(\mathcal{A}_{12}\right) \).

\textbf{Claim 4}: \( \phi\left(\mathcal{A}_{21}\right) \subseteq \alpha\left(\mathcal{A}_{21}\right)\) .

Since \( \phi\left(P_{1}\right) \subseteq \alpha\left(\mathcal{A}_{11}\right) \) and  \(\alpha\left(P_{1}\right) \phi\left(P_{1}\right) \alpha\left(P_{1}\right)=0 \) according  to Equation 11  yields  \(\phi\left(P_{1}\right)=0  \). We consider  \(A_{11}=P_{1} \) and \(  A_{2 i}=A_{21} \) in Equation 3, associated with \( \phi\left(\mathcal{A}_{11}\right) \subseteq   \alpha\left(\mathcal{A}_{11}\right) \) and  \(\phi\left(P_{1}\right)=0 \), we obtain that
\begin{equation}
\phi\left(A_{12} A_{21}\right)=\phi\left(A_{12}\right) \alpha\left(A_{21}\right)+\alpha\left(A_{12}\right) \phi\left(A_{21}\right) .
\end{equation}

We multiply the right-hand side of Equation 20 by \( Q_{2}=\alpha\left(P_{2}\right)  \)and yield \( \alpha\left(A_{12}\right) \phi\left(A_{21}\right) \alpha\left(P_{2}\right)=0 \),   which   in   turn   gives \(\alpha\left(P_{2}\right) \phi\left(A_{21}\right) \alpha\left(P_{2}\right)=0 \). \(\alpha\left(P_{1}\right) \phi\left(A_{21}\right)=0 \) and Equation 7, implies that  \(\phi\left(\mathcal{A}_{21}\right) \subseteq \alpha\left(\mathcal{A}_{21}\right)\) .

Equation 5 obviously holds when \( j \neq s \) according to Claims 1-4. Next, it is only necessary to prove the case when \( j=s \), which is the following eight cases:
\begin{equation*}
(1)  \phi\left(A_{11} B_{11}\right)=\phi\left(A_{11}\right) \alpha\left(B_{11}\right)+\alpha\left(A_{11}\right) \phi\left(B_{11}\right) ;
\end{equation*}
\begin{equation*}
(2)  \phi\left(A_{11} B_{12}\right)=\phi\left(A_{11}\right) \alpha\left(B_{12}\right)+\alpha\left(A_{11}\right) \phi\left(B_{12}\right) ;
\end{equation*}
\begin{equation*}
(3)  \phi\left(A_{12} B_{21}\right)=\phi\left(A_{12}\right) \alpha\left(B_{21}\right)+\alpha\left(A_{12}\right) \phi\left(B_{21}\right) ;
\end{equation*}
\begin{equation*}
(4)  \phi\left(A_{12} B_{22}\right)=\phi\left(A_{12}\right) \alpha\left(B_{22}\right)+\alpha\left(A_{12}\right) \phi\left(B_{22}\right) ;
\end{equation*}
\begin{equation*}
(5)  \phi\left(A_{22} B_{21}\right)=\phi\left(A_{22}\right) \alpha\left(B_{21}\right)+\alpha\left(A_{22}\right) \phi\left(B_{21}\right) ;
\end{equation*}
\begin{equation*}
(6)  \phi\left(A_{22} B_{22}\right)=\phi\left(A_{22}\right) \alpha\left(B_{22}\right)+\alpha\left(A_{22}\right) \phi\left(B_{22}\right) ;
\end{equation*}
\begin{equation*}
(7)  \phi\left(A_{21} B_{11}\right)=\phi\left(A_{21}\right) \alpha\left(B_{11}\right)+\alpha\left(A_{21}\right) \phi\left(B_{11}\right) ;
\end{equation*}
\begin{equation*}
(8)  \phi\left(A_{21} B_{12}\right)=\phi\left(A_{21}\right) \alpha\left(B_{12}\right)+\alpha\left(A_{21}\right) \phi\left(B_{12}\right) .
\end{equation*}
Case (1), case (2), and case (3) can be obtained from Equation 18, Equation 17, and  Equation 20, respectively.

Since \( B_{22}=\mu Q_{2} \) for  some \( \mu \in \mathbb{F} \) and  \(\phi\left(B_{22}\right)=0\), then  case   (4) holds. Similarly, it follows that case (5) holds. Additionally, case (6) is trivial because  \(\phi\left(\mathcal{A}_{22}\right)=\{0\}\).

According to case (1) and case (2) we can get:

\begin{equation*}
\begin{aligned}
\phi\left(A_{12} A_{21} B_{11}\right) & =\phi\left(A_{12}\right) \alpha\left(A_{21} B_{11}\right)+\alpha\left(A_{12}\right) \phi\left(A_{21} B_{11}\right) \\
&=\phi\left(A_{12} A_{21}\right) \alpha\left(B_{11}\right)+\alpha\left(A_{12} A_{21}\right) \phi\left(B_{11}\right) \\
&=\phi\left(A_{12}\right) \alpha\left(A_{21} B_{11}\right)+\alpha\left(A_{12}\right) \phi\left(A_{21}\right) \alpha\left(B_{11}\right)+\alpha\left(A_{12} A_{21}\right) \phi\left(B_{11}\right)
\end{aligned}
\end{equation*}

Claim 4 implies case (7).

Case (8) can be derived similarly from the following equation:

\begin{equation*}
\begin{aligned}
\phi\left(A_{12} A_{21} B_{12}\right) & =\phi\left(A_{12}\right) \alpha\left(A_{21} B_{12}\right)+\alpha\left(A_{12}\right) \phi\left(A_{21} B_{12}\right) \\
&=\phi\left(A_{12} A_{21}\right) \alpha\left(B_{12}\right)+\alpha\left(A_{12} A_{21}\right) \phi\left(B_{12}\right) \\
&=\phi\left(A_{12}\right) \alpha\left(A_{21} B_{12}\right)+\alpha\left(A_{12}\right) \phi\left(A_{21}\right) \alpha\left(B_{12}\right)+\alpha\left(A_{12} A_{21}\right) \phi\left(B_{12}\right).
\end{aligned}
\end{equation*}

So \( \phi \) is an  \((\alpha, \alpha)\)-derivation on  \(\mathcal{A}\), and \( \delta \) is an  \((\alpha, \alpha)\)-derivation on \( \mathcal{A} \).

\textbf{Corollary 2.1.} Suppose that \( G \) is a nonzero element in  \(\mathcal{A}=B(X)\), \(\alpha \) is an automorphism on \( \mathcal{A}\), and  \(\delta \) is a linear mapping on  \(\mathcal{A} \). If \( \delta \) is \( (\alpha, \alpha) \)-generalized derivable at \( G \), which means that \( \delta(G)=\delta(A) \alpha(B)+\alpha(A) \delta(B)-\alpha(A) \delta(I) \alpha(B) \) for all \( A, B \in \mathcal{A}  \) with \( A B=G \), then  \(\delta\)  is an \( (\alpha, \alpha) \)-generalized derivation.

\textbf{Proof.} We define a new linear mapping \(\phi: \mathcal{A} \rightarrow \mathcal{A} \) by  \(\phi(A)=\delta(A)-\delta(I) \alpha(A)\). From the equations
\begin{equation*}
\phi(A B)=\delta(A B)-\delta(I) \alpha(A B)
\end{equation*}
and
\begin{equation*}
\phi(A) \alpha(B)+\alpha(A) \phi(B)=(\delta(A)-\delta(I) \alpha(A)) \alpha(B)+\alpha(A)(\delta(B)-\delta(I) \alpha(B)),
\end{equation*}

 \(\delta \) is \( (\alpha, \alpha) \)-generalized derivable at \( G \) if and only if \( \phi \) is \( (\alpha, \alpha)\) -derivable at \( G \). Then, \( \delta \) is \( (\alpha, \alpha)\)-generalized derivable if and only if \( \phi  \) is \( (\alpha, \alpha)\)-derivable.

Because \( \phi \) is  \((\alpha, \alpha) \)-derivable at \( G\) , by Theorem 2.1, \( \phi \) is  \((\alpha, \alpha) \)-derivable. Then \( \delta \) is an \( (\alpha, \alpha)\) -generalized derivation on  \(\mathcal{A}\) .

\bibliographystyle{unsrt}  


\end{document}